\documentclass[a4paper]{amsart}

\usepackage{amsmath}
\usepackage{amsthm}
\usepackage{amssymb}
\usepackage{graphicx}
\usepackage{url}

\newcommand{\R}{\mathbb{R}}
\newcommand{\N}{\mathbb{N}}

\newcommand{\A}{\mathcal{A}}

\newcommand{\dt}{\Delta t}

\newcommand{\pa}{\partial}
\newcommand{\dom}{\mathcal{D}}
\newcommand{\f}[2]{\frac{#1}{#2}}
\newcommand{\abs}[1]{\left\lvert#1\right\rvert}
\newcommand{\norm}[1]{\left\lVert#1\right\rVert}

\newcommand{\floor}[1]{\left\lfloor#1\right\rfloor}

\newcommand{\ds}{\displaystyle}
\newcommand{\s}{\mathfrak{S}}

\DeclareMathOperator{\supp}{supp}
\DeclareMathOperator{\esupp}{esupp}

\theoremstyle{plain}
\newtheorem{theorem}{Theorem}[section]
\newtheorem{proposition}[theorem]{Proposition}

\theoremstyle{definition}
\newtheorem{definition}[theorem]{Definition}
\newtheorem{example}[theorem]{Example}
\newtheorem{assumption}[theorem]{Assumption}
\newtheorem{remark}[theorem]{Remark}

\begin{document}

\title{Cubature on Wiener space in infinite dimension}

\author{Christian Bayer and Josef Teichmann}

\address{Vienna University of Technology, Department of mathematical methods in
  economics, Research Group e105 Financial and Actuarial Mathematics,
  Wiedner Hauptstrasse 8-10, A-1040 Wien, Austria}
\email{cbayer@fam.tuwien.ac.at, jteichma@fam.tuwien.ac.at}
\thanks{The authors gratefully acknowledge the support from the FWF-grant Y328}

\begin{abstract}
  We prove a stochastic Taylor expansion for SPDEs and apply this result to
  obtain cubature methods, i.e.~high order weak approximation schemes for
  SPDEs, in the spirit of T.~Lyons and N.~Victoir. We can prove a high-order
  weak convergence for well-defined classes of test functions if the process
  starts at sufficiently regular initial values. We can also derive analogous
  results in the presence of L\'evy processes of finite type, here the results
  seem to be new even in finite dimension. Several numerical examples are
  added.
\end{abstract}
\maketitle

\section{Introduction}
\label{sec:introduction}

Let $H$ be a real separable Hilbert space. We consider the stochastic
(partial) differential equation for a diffusion process $ X $ with values in
$H$
\begin{equation}
  \label{eq:sde}
  dX_t^x = (AX^x_{t} + \alpha(X^x_{t})) dt + \sum_{i=1}^d
  \beta_i(X^x_{t}) dB^i_t
\end{equation}
or -- in the presence of jumps -- the stochastic differential equation for a
jump-diffusion process $X$ with values in $H$
\begin{equation}
  \label{eq:sde-j}
  dX_t^x = (AX^x_{t^-} + \alpha(X^x_{t^-})) dt + \sum_{i=1}^d
  \beta_i(X^x_{t^-}) dB^i_t +
  \sum_{j=1}^{e}\delta_{j}(X_{t^-}^{x})dL_{t}^{j}, 
\end{equation}
where, in general, $A:\mathcal{D}(A)\subset H \to H$ denotes an unbounded
linear operator, $\alpha, \beta_1,\ldots,\beta_d,\delta_1,\ldots,\delta_e:H
\to H$ denote $C^\infty$-bounded vector fields -- i.e.~the vector fields are
smooth and all the derivatives (of degree $\geq1$) are bounded -- and
$(B_t)_{t \geq 0} = (B^1_t,\ldots,B^d_t)_{t \geq 0}$ denotes a finite
dimensional Brownian motion on the Wiener space $(\Omega,\mathcal{F},P)$, and
$(L^j_t)_{t \geq0 }$ a compound Poisson process with jump-rate $ \mu_j$ for
$j=1,\ldots,e$. The initial value $x\in H$ appears as superscript in the
notation $X^x_t$ of the solution process, i.e.~$X^x_0 = x$. We assume that $A$
is the generator of a $C_0$-semigroup denoted by $(S_t)_{t \geq 0}$. The main
reference for equations of the type~\eqref{eq:sde} is the monograph of G.~da
Prato and J.~Zabczyk~\cite{dap/zab92}. In the (general) L\'{e}vy case we refer
to \cite{fil/tap07} and \cite{peszab:07}, even though we do not need results
of their strength in our case since all L\'evy processes under consideration
are of finite type.

Due to the unboundedness of the operator $A$ several concepts of solutions
to~\eqref{eq:sde-j} arise. We refer for the precise definitions to the
excellent monograph \cite{peszab:07}.  The most direct analogue to the finite
dimensional setting is the concept of a \emph{strong solution}, which is
defined by
\begin{multline}
  \label{eq:strong_solution}
  X^x_t = x + \int_0^t (AX^x_{s^-} + \alpha(X^x_{s^-})) ds + \sum_{i=1}^d
  \int_0^t \beta_i(X^x_{s^-}) dB^i_s + \\ +\sum_{j=1}^{e} \int_0^t
  \delta_{j}(X_{s^-}^{x})dL_{s}^{j},
\end{multline}
i.e.~by the integrated version of equation~\eqref{eq:sde}. Besides the obvious
integrability assumptions we need in particular that $X^x_t\in\mathcal{D}(A)$
for all $t\geq0$ almost surely with respect to $P$. More relevant is the
concept of a \emph{mild solution}, which is related to \eqref{eq:sde-j} via
variation of constants. Indeed, a mild solution is a process $X^x_t$
satisfying
\begin{multline}
  \label{eq:mild_solution}
  X^x_t = S_t x + \int_0^t S_{t-s} \alpha(X_{s^-}^x) ds + \sum_{i=1}^d
  \int_0^t
  S_{t-s} \beta_i(X^x_{s^-}) dB^i_s + \\
  + \sum_{j=1}^{e}\int_0^t S_{t-s} \delta_{j}(X_{s^-}^{x})dL_{s}^{j}
\end{multline}
given the obvious integrability assumptions. By It\^{o}'s formula, any strong
solution is a mild solution, but not the other way round. Note that mild
solutions need not be semi-martingales, because there is no semi-martingale
decomposition if the process evolves outside of
$\mathcal{D}(A)$. Consequently, a Stratonovich formulation of~\eqref{eq:sde-j}
is, in general, not possible for mild solutions.

Of course, neither strong nor mild solutions can usually be given explicitly,
which makes numerical approximation necessary. We are interested in \emph{weak
  approximation} of the solution in the sense that we want to approximate the
value
\begin{equation*}
  P_t f(x) = E(f(X^x_t))
\end{equation*}
for a suitable class of test functions $f: H \to \R$ at initial values $ x \in
H $. It is well-known that the function $(t,x) \mapsto P_tf(x)$ solves the
Kolmogorov equation in the weak sense, see for instance \cite{dap/zab92} in
the diffusion case.

Let us assume for a moment that there are no jump-components: usually,
infinite dimensional SDEs are numerically solved by finite element or finite
difference sche\-mes, see, for instance, \cite{hau03}, \cite{yan05},
\cite{gyo98} and \cite{gyo99}. For HJM models of financial mathematics a
finite difference scheme and a finite element scheme have been implemented in
\cite{bjo/sze/tem/zou08}. This means that the original equation is projected
onto some finite dimensional subspace $H_h\subset H$ and $A$ is approximated
by some operator $A_h$ defined thereon. This procedure, which corresponds to a
space discretization of the stochastic PDE, is followed by a conversion of the
stochastic differential equation on $H_h$ to a stochastic difference equation
by discretizing in time, using an Euler method or a related scheme. Finally,
the stochastic difference equation is solved by Monte-Carlo simulation, which
may be interpreted as a discretization on the Wiener space. For general
information about approximation of finite dimensional SDEs see
\cite{kloe/pla92}.

We want to tackle the problem in the reverse order: we want to do the
discretization on the Wiener space $\Omega$ first, reducing the problem to a
deterministic problem, i.e., one replaces the $d$-dimensional Brownian motion
with finitely many trajectories of bounded variation chosen with well-defined
probabilities such that certain moments (of iterated Stratonovich integrals)
match.  The resulting deterministic problem can be solved by standard methods
for the numerical treatment of deterministic PDEs, e.g.~by standard finite
element or finite difference methods. The benefit is that once the
discretization on the Wiener space has been done, we can use the
well-established theory of the corresponding deterministic PDE-problems,
without any complications from stochasticity. Our method of choice for
discretization on $\Omega$ is ``\emph{cubature on Wiener space}'', developed
by Terry Lyons and Nicolas Victoir in \cite{lyo/vic04} and by Shigeo Kusuoka
in \cite{kus01} and \cite{kus04}, see also \cite{tei06}.  In the spirit of
these methods we shall obtain weak approximation schemes of any prescribed
order of convergence. Notice here that we discretize in the presence of the
unbounded operator $A$ in the drift vector field. Certainly our Assumptions
\ref{ass:dainfty} seem very restrictive, but these assumptions are the
appropriate analogues of the assumptions in finite dimension that the vector
fields are bounded, $C^{\infty}$-bounded (see \cite{kus04} and
\cite{lyo/vic04}).

Before going into details, let us motivate the use of cubature formulas in the
present context. Let $X^x_t(\omega)$ denote the solution of~\eqref{eq:sde},
formally rewritten in Stratonovich form, if each ``$dB^i_t$'' is replaced by
``$d\omega^i(t)$'', i.e.~
\begin{multline}
  \label{eq:cubature_ode}
  dX^x_t(\omega) = (AX^x_t(\omega) + \alpha(X^x_t(\omega)) - \frac{1}{2}
  \sum_{i=1}^d D\beta_i(X_t^x(\omega)) \cdot \beta_i (X_t^x(\omega))) dt + \\
  + \sum_{i=1}^d \beta_i(X^x_t(\omega)) d\omega^i(t)
\end{multline}
for a curve function $\omega = (\omega^1,\ldots,\omega^d) : \R_{\geq0} \to
\R^d$ of bounded variation. Roughly speaking the idea of cubature on Wiener
space is to construct short-time asymptotics (for some given degree of
accuracy $ m \geq 2$)
\begin{equation*}
  E(f(X^x_t)) = P_tf(x) = \sum_{l=1}^N \lambda_l f(X^x_t(\omega_l)) +
  \mathcal{O}(t^{\frac{m+1}{2}}), 
\end{equation*}
with some positive, time-independent weights $\lambda_1,\ldots,\lambda_N$
satisfying $\lambda_1+\cdots+\lambda_N=1$ and some $d$-dimensional paths
$\omega_1,\ldots, \omega_N$ of bounded variation. Of course, the weights and
paths are chosen in a specific way, which will be described later in more
detail, and the asymptotics will only hold for some class of test functions
$f$. Notice in particular that the cubature paths $\omega_1,\ldots,\omega_N$
depend on the interval $[0,t]$ -- they become rougher as $t$ approaches
$0$. The aforementioned procedure replaces the SDE by $N$ deterministic,
well-defined PDEs, which have unique mild solutions. The iteration of the
short-time asymptotics due to the Markov property then yields a weak, high
order approximation scheme.

\emph{Here also the main advantage of cubature methods in contrast to Taylor
  methods gets visible}. The time-discretization in the realm of cubature
methods always leads to reasonable expressions, namely to \emph{reasonable
  partial differential equations of type~\eqref{eq:cubature_ode}}. If we
wanted to apply the usual discretization methods in time like the
Euler-Maruyama method, we might run into problems. Indeed, the naive Euler
scheme is well-suited for the differential formulation~\eqref{eq:sde} of the
problem,
\begin{equation*}
  X_0 = x \text{ and } X_n = (AX_{n-1} + \alpha(X_{n-1}))\frac{t}{n} +
  \sum_{i=1}^d \beta_i(X_{n-1}) \Delta_n B^i, 
\end{equation*}
for $ n \geq 1 $, but it might immediately lead to some $
X_n\notin\mathcal{D}(A)$. Even in the case of a well-defined strong solution,
there is no reason why the discrete approximation should always stay in
$\mathcal{D}(A)$. Hence the naive implementation of the Euler-Maruyama method
does not work.

Only formulation~\eqref{eq:mild_solution} seems to be suitable for using an
Euler-like method. If one understands the semigroup $ S $ well, one can
approximate $X^x_t$ by expressing the integrals in~\eqref{eq:mild_solution} as
Riemannian sums, involving evaluations of $S_{t-s}$, which yields a sort of
strong Euler method (see for instance the very interesting book
\cite{preroe07} for strong convergence theorems in this direction).

We do not discretize the integral in formulation~\eqref{eq:mild_solution}, but
(weakly) approximate the Brownian motion's paths by paths of bounded variation
such that we obtain a weak approximations of $X^x_t$.

The presence of jumps does not lead to more complicated expressions, since the
short time asymptotics of a jump-diffusion can be easily derived from a
diffusion's short-time asymptotics by conditioning on the jumps. The arising
picture is the following. Discretizing the equation~\eqref{eq:sde-j} means to
allow a certain number of jumps between to consecutive points in the time
grid. Between two jumps we apply a diffusion cubature formula to express the
short-time asymptotics. This yields as a corollary of the diffusion theory
also the jump diffusion theory.

\begin{remark}
  \label{rem:infinite-activity}
  A direct application of the cubature on Wiener space technique to jump
  diffusions driven by L\'evy processes of infinite activity is not
  possible. However, notice that any L\'evy process can be approximated by
  processes with finite activity in a weak sense, see, for instance,
  \cite{contan:04} and \cite{peszab:07}, and therefore the solutions of the
  corresponding Hilbert-space valued SDEs converge weakly. In this sense,
  cubature methods can be use for jump diffusions driven by L\'evy processes
  of infinite activity, too.
\end{remark}

The article is organized as follows. In Section \ref{sec:setting} we describe
the analytic setting for a stochastic Taylor expansion to work. This is a
delicate question since we deal with one unbounded vector field. In Section
\ref{sec:cubat-bound-oper} and Section \ref{sec:cubat-equat-mathc} we work out
the cubature method from the scratch and prove the relevant convergence
results in the diffusion case. In Section \ref{sec:jump-structure} we allow
for jumps and prove the associated short-time asymptotics which is relevant to
set up a weak approximation scheme. In Section \ref{sec:numerical-examples} we
apply our method to several examples to demonstrate the results.

\section{Setting and assumptions}
\label{sec:setting}

Let $(\Omega,\mathcal{F},P,(\mathcal{F}_{t})_{t\geq0})$ be a filtered
probability space with a filtration $(\mathcal{F}_{t})_{t\geq0}$ satisfying
the usual conditions. Let $(B_{t})_{t\geq0}$ be a $d$-dimensional Brownian
motion and $(L_{t}^{j})_{t\geq0}$, $j = 1, \ldots, e$ be $e$ independent
compound Poisson processes given by
\begin{equation*}
  L_t^j := \sum_{k = 1}^{N_t^j} Z_k^j,\label{poisson-process}
\end{equation*}
where $N_t^j$ denotes a Poisson process with jump rate $\mu_{j}>0$ and $Z^j =
(Z_k^j)_{k \geq 1}$ is an i.i.d.~sequence of random variables with
distribution $\nu_{j}$ for $j=1,\dots,e$, such that each $\nu_{j}$ admits all
moments. All notions are with respect to the given filtration. We assume
furthermore that all sources of randomness are mutually independent.

Let $ H $ be a separable Hilbert space. We furthermore fix a strongly
continuous semi-group $ S $ on $H$ with generator $ A $.  Let
$\alpha,\beta_{1},\dots,\beta_{d}$, the diffusion vector fields, and $\delta
_{1},\dots,\delta_{e}$, the jump vector fields, be $C^{\infty}$-bounded on
$H$, that is, the vector fields are infinitely often differentiable with
bounded partial derivatives of all proper orders $ n \geq 1 $. We consider the
mild c\`{a}dl\`{a}g solution $(X_{t}^{x})_{t \geq 0}$ of a stochastic
differential equation
\begin{align}
  dX_{t}^{x} & =(AX^x_{t^-} + \alpha(X_{t^-}^{x}))dt +
  \sum_{i=1}^{d}\beta_{i}(X_{t^-}^{x})dB_{t}^{i} +
  \sum_{j=1}^{e}\delta_{j}(X_{t^-}^{x})dL_{t}^{j},\\
  X_{0}^{x} & =x \in H.
\end{align}
See \cite{fil/tap07} for all necessary details on existence and uniqueness for
the previous equation. Notice furthermore the decomposition theorem in
\cite{for/lue/tei07}, which states that we do not need any further existence
and uniqueness results in this case: in particular, we do not need to impose
further (contractivity) conditions on $A$ as in~\cite{fil/tap07} in the finite
activity case.
 
The previous conditions are slightly more than standard for existence and
uniqueness of mild solutions, i.e.~in \cite{fil/tap07} the authors need
Lipschitz conditions on the vector fields, whereas we assume them to be
$C^{\infty}$-bounded. In order to formulate a stochastic Taylor expansion we
shall need one main assumption, which we formulate in the sequel. This
assumption has already been successfully applied in several circumstances,
e.g.~\cite{bau/tei05}, \cite{fil/tei03} or the recent \cite{for/lue/tei07}.

We apply the following notations for Hilbert spaces $ \dom(A^k) $,
\begin{equation*}
  \begin{array}{ccl}
    \dom(A^0) &:=& H,\\
    \dom(A^k)&:=&\ds \left\{h\in H| \, h\in \dom(A^{k-1})
      \hspace{2mm}\mbox{and}\hspace{2mm}
      A^{k-1}h\in \dom(A)\right\},\\
    ||h||^2_{\dom(A^k)}&:=&\ds \sum_{i=0}^k ||A^ih||^2,\\
    \dom(A^{\infty})&=&\ds \bigcap_{k\geq 0} \dom(A^k),
  \end{array}
\end{equation*}
which we need in order to specify the main analytic condition for our
considerations:

\begin{assumption}\label{ass:dainfty}
  We assume that $\alpha,\beta_{1},\dots,\beta_{d}$, the diffusion vector
  fields, and $\delta_{1},\dots,\delta_{e}$, the jump vector fields, map
  $\dom(A^k) \to \dom(A^k)$ and are $C^{\infty}$-bounded thereon for each $ k
  \geq 0 $, that is, the vector fields are infinitely often differentiable
  with bounded partial derivatives of all orders $n \geq 1 $ on the Hilbert
  space $\dom(A^k)$ for each $ k \geq 0 $.
\end{assumption}

\begin{remark}
  These assumptions are the appropriate analogues of the assumptions in finite
  dimension that the vector fields are bounded, $C^{\infty}$-bounded (see
  \cite{kus04} and \cite{lyo/vic04}). In order to establish a true convergence
  rate one needs an additional cut-off argument, which is outlined in Remark
  \ref{remark_ass:dainf}. This can -- like in the finite dimensional setting
  -- certainly be improved.
\end{remark}

\begin{example}\label{banach-map-vector-fields}
  In order to show examples of vector fields which are $ C^{\infty}$-bounded
  on $ \dom(A^k) $ consider the following structure. Let $H$ be a separable
  Hilbert space and $ A $ the generator of a strongly continuous semigroup. We
  know that $ \dom(A^{\infty}) $ is a Fr\'echet space and an injective limit
  of the Hilbert spaces $\dom(A^k) $ for $ k \geq 0$. Following the analysis
  as developed in \cite{fil/tei03} (see also \cite{KriMic97} and \cite{Ham82}
  were the analytic concepts have been originally developed), we can consider
  the vector field $ V: U \subset H \to \dom(A^{\infty}) $. If $V$ is smooth
  in the sense explained in \cite{fil/tei03} and has the property that its
  derivatives of order $ n \geq 1 $ are bounded on $ U \subset H $, then $ V$
  is obviously a $ C^{\infty}$-bounded vector field and additionally $
  V|_{\dom(A^{\infty})} $ is a Banach-map-vector field in the sense of
  \cite{fil/tei03}. Such vector fields constitute a class, where the above
  assumptions can be readily checked.
\end{example}

\section{The case when $A$ is bounded linear}
\label{sec:cubat-bound-oper}

We shall assume in this section that there are no jumps, i.e.~we consider
\begin{equation}
  \label{eq:sde_nojump}
  dX^x_t = (AX^x_t + \alpha(X^x_t)) dt + \sum_{i=1}^d \beta_i(X^x_t) dB^i_t.
\end{equation} 
In order to give an introduction to cubature on Wiener space, we consider the
problem for a bounded operator $A$. In this case, there are virtually no
differences to the finite dimensional setting, except the fact that the drift
vector field does not need to be bounded by some constant on the whole Hilbert
space (due to the presence of one linear operator in it). Remember that in
\cite{lyo/vic04} and \cite{kus01} and \cite{kus04} one deals with globally
bounded vector fields. We shall circumvent this problem by a small refinement
of the arguments.

Since mild and strong solutions coincide, we can always work with strong
solutions, which are semi-martingales. Consequently, we can
rewrite~\eqref{eq:sde_nojump} into its Stratonovich form
\begin{equation}
  \label{eq:sde_stratonovich}
  dX^x_t = \beta_0(X^x_t) dt + \sum_{i=1}^d \beta_i(X^x_t) \circ dB^i_t,
\end{equation}
where $\beta_0: H\to H$ denotes the Stratonovich-corrected drift, i.e.~
\begin{equation}
  \label{eq:stratonovich_drift}
  \beta_0(x) = Ax + \alpha(x) - \f{1}{2}\sum_{i=1}^d D\beta_i(x)\cdot
  \beta_i(x), 
\end{equation}
where
\begin{equation*}
  DF(x) \cdot v = \left. \f{\pa}{\pa \epsilon}\right|_{\epsilon=0}
  F(x+\epsilon v)
\end{equation*}
denotes the Fr\'{e}chet derivative of a function or vector field $F$.  This
notation enables us to write
\begin{equation*}
  dX^x_t = \sum_{i=0}^d \beta_i(X^x_t)\circ dB^i_t,
\end{equation*}
where we use the convention that ``$\circ dB^0_t=dt$''.

The following notions form the core of cubature on Wiener space. We only give
a short description and refer the reader to \cite{lyo/vic04} and \cite{tei06}
for more details.  Let $\mathcal{A}$ denote the set of all multi-indices in
$\{0,\ldots,d\}$. We define a degree on $\mathcal{A}$ by setting
\begin{equation*}
  \deg(i_1,\ldots,i_k) = k + \#\{1\leq j\leq k\, | \, i_j = 0\},
\end{equation*}
$k\in\N$, $(i_1,\ldots,i_k)\in\mathcal{A}$. We have to count all the zeros
twice because of the different scalings for $t=B^0_t$ and the Brownian motion.

Recall that any vector field $\beta$ can be interpreted as a first order
differential operator on test functions $f$ by
\begin{equation*}
  (\beta f)(x) = Df(x)\cdot\beta(x),\quad x\in H.
\end{equation*}
For a multi-index $(i_1,\ldots,i_k)\in\mathcal{A}$, $k\in\N$, let
\begin{equation*}
  B^{(i_1,\ldots,i_k)}_t = \int_{0 \leq t_1 \leq \cdots \leq t_k \leq t} \circ
  dB^{i_1}_{t_1}\cdots \circ dB^{i_k}_{t_k}
\end{equation*}
denote the corresponding iterated Stratonovich integral. The iterated
integrals form the building blocks of the stochastic Taylor formula, see
\cite{bau04}.
\begin{proposition}[Stochastic Taylor expansion]
  \label{prop:stoch_taylor}
  Let $f\in C^\infty(H;\R)$ and fix $0 < t < 1$, $m\in\N$, $m\geq1$, $x\in
  H$. Then we have
  \begin{equation*}
    f(X^x_t) = \sum_{\substack{(i_1,\ldots,i_k)\in\mathcal{A}\\
        \deg(i_1,\ldots,i_k) \leq m}} (\beta_{i_1}\cdots \beta_{i_k}f)(x)
    B^{(i_1,\ldots,i_k)}_t + R_m(t,f,x),
  \end{equation*}
  with
  \begin{equation*}
    \sup_{x\in H} \sqrt{E(R_m(t,f,x)^2)} \leq C
    t^{\f{m+1}{2}} \max_{m<\deg(i_1,\ldots,i_k)\leq
      m+2}\; \sup_{y \in H}|\beta_{i_1}\cdots\beta_{i_k}f(y)|.
  \end{equation*}
\end{proposition}
\begin{remark}\label{remark_remainder_finite}
  Notice that under the assumptions of Proposition \ref{prop:stoch_taylor} the
  bound for the remainder term can be infinite.  However, if $ f $ is bounded,
  $C^{\infty}$-bounded and has bounded support, then we can guarantee that
$$
\sup_{y \in H}|\beta_{i_1}\cdots\beta_{i_k}f(y)| < \infty.
$$
In the unbounded case this question is more subtle, see Example
\ref{example_ass:dainf} and Remark \ref{remark_ass:dainf}.
\end{remark}
Proposition~\ref{prop:stoch_taylor} shows that iterated Stratonovich integrals
play the r\^{o}le of polynomials in the stochastic Taylor expansion.
Consequently, it is natural to use them in order to define cubature
formulas. Let $C_{bv}([0,t];\R^d)$ denote the space of paths of bounded
variation taking values in $\R^d$. As for the Brownian motion, we append a
component $\omega^0(t)=t$ for any $\omega\in C_{bv}([0,t];\R^d)$. Furthermore,
we establish the following convention: whenever $X^x_t$ is the solution to
some stochastic differential equation of type~(\ref{eq:sde}) driven by $B$,
whether on a finite or infinite dimensional space, and $\omega\in
C_{bv}([0,t];\R^d)$, we denote by $X^x_t(\omega)$ the solution of the
deterministic differential equation given by formally replacing all
occurrences of ``$\circ dB^i_s$'' with ``$d\omega^i(s)$'' (with the same
initial values), see equation~(\ref{eq:cubature_ode}) as compared
to~(\ref{eq:sde}). Note that it is necessary that the SDE for $X$ is formally
formulated in the Stratonovich sense (recall that the Stratonovich formulation
does not necessarily make sense for \eqref{eq:sde}).
\begin{definition}
  \label{def:cubature_wiener_space}
  Fix $t>0$ and $m\geq1$. Positive weights $\lambda_1,\ldots,\lambda_N$
  summing up to $1$ and paths $\omega_1,\ldots,\omega_N\in C_{bv}([0,t];\R^d)$
  form a \emph{cubature formula on Wiener space} of degree $m$ if for all
  multi-indices $(i_1,\ldots,i_k)\in \mathcal{A}$ with
  $\deg(i_1,\ldots,i_k)\leq m$, $k\in\N$, we have that
  \begin{equation*}
    E(B^{(i_1,\ldots,i_k)}_t) = \sum_{l=1}^N \lambda_l
    B^{(i_1,\ldots,i_k)}_t(\omega_l), 
  \end{equation*}
  where we used the convention in line with the previous one, namely
$$
B^{(i_1,\ldots,i_k)}_t(\omega)= \int_{0 \leq t_1 \leq \cdots \leq t_k \leq t}
d \omega^{i_1}(t_1) \cdots d \omega^{i_k}(t_k).
$$ 
\end{definition}

Lyons and Victoir \cite{lyo/vic04} show the existence of cubature formulas on
Wiener space for any $d$ and size $N \leq \#\{I\in\mathcal{A} | \deg(I)\leq
m\}$ by applying Chakalov's theorem on cubature formulas and Chow's theorem
for nilpotent Lie groups. Moreover, due to the scaling properties of Brownian
motion (and its iterated Stratonovich integrals), i.e.
\begin{equation*}
  B^{(i_1,\ldots,i_k)}_t =^{\text{law}} \sqrt{t}^{\deg(i_1,\ldots,i_k)}
  B^{(i_1,\ldots,i_k)}_1, 
\end{equation*}
it is sufficient to construct cubature paths for $t=1$.

\begin{assumption}
  \label{ass:one-cubature-formula}
  Once and for all, we fix one cubature formula $\widetilde{\omega}_1, \ldots,
  \widetilde{\omega}_N$ with weights $\lambda_1, \ldots, \lambda_N$ of degree
  $m \geq 1$ on the interval $[0,1]$. By abuse of notation, for any $t>0$, we
  will denote $\omega_l(s) = \widetilde{\omega}_l(s/\sqrt{t})$, $s \in [0,t]$,
  $l = 1,\ldots, N$, the corresponding cubature formula for $[0,t]$.
\end{assumption}

\begin{example}\label{ex:cubature_paths}
  For $d=1$ Brownian motions, a cubature formula on Wiener space of degree
  $m=3$ is given by $N=2$ paths
  \begin{equation*}
    \omega_1(s) = -\f{s}{\sqrt{t}},\ \omega_2(s) = \f{s}{\sqrt{t}}
  \end{equation*}
  for fixed time horizon $t$. The corresponding weights are given by
  $\lambda_1 = \lambda_2 = \f{1}{2}$.
\end{example}

Combining the stochastic Taylor expansion, the deterministic Taylor expansion
for solutions of ODEs on $H$ and a cubature formula on Wiener space yields a
one-step scheme for weak approximation of the SDE~(\ref{eq:sde}) for bounded
operators $A$ in precisely the same way as in \cite{lyo/vic04}. Indeed, we get
\begin{multline}
  \label{eq:bounded_1step}
  \sup_{x\in H}\abs{E(f(X^x_t)) - \sum\nolimits_{l=1}^N \lambda_l
    f(X^x_t(\omega_l))} \\
  \leq C t^{\f{m+1}{2}} \max_{\substack{ (i_1,\ldots,i_k)\in \mathcal{A} \\ m
      < \deg(i_1,\ldots,i_k)\leq m+2}} \sup_{y \in
    H}|\beta_{i_1}\cdots\beta_{i_k}f(y)|,
\end{multline}
for $0<t<1$, $f\in C^{\infty}(H)$.

For the multi-step method, divide the interval $[0,T]$ into $p$ subintervals
according to the partition $0= t_0 < t_1 < \cdots < t_p = T$. For a
multi-index $(l_1,\ldots,l_p)\in \{1,\ldots,N\}^p$ consider the path
$\omega_{l_1,\ldots,l_p}$ defined by concatenating the paths
$\omega_{l_1},\ldots,\omega_{l_p}$, i.e.~$\omega_{l_1,\ldots,l_p}(t) = \omega_
{l_1}(t)$ for $t \in [0, t_1[$ and
\begin{equation*}
  \omega_{l_1,\ldots,l_p}(t) = \omega_{l_1,\ldots,l_p}(t_{r-1}) + 
  \omega_{l_r}(t-t_{r-1})
\end{equation*}
for $r$ such that $t\in [t_{r-1},t_r[$, where $\omega_{l_r}$ is scaled to be a
cubature path on the interval $[0,t_r-t_{r-1}]$.
\begin{proposition}
  \label{prop:bounded_multistep}
  Fix $T>0$, $m\in \N$, a cubature formula of degree $m$ as in
  Definition~\ref{def:cubature_wiener_space} and a partition of $[0,T]$ as
  above. For every $ f \in C^{\infty}(H) $ with
  \begin{equation*}
    \sup_{0\leq t \leq T} \sup_{y \in H}
    \abs{\beta_{i_1}\cdots \beta_{i_k} (P_tf)(y)} < \infty
  \end{equation*}
  for all $(i_1,\ldots,i_k)\in\mathcal{A}$ with $m< \deg(i_1,\ldots,i_k) \leq
  m+2$, $k\in\N$, there is a constant $D$ independent of the partition such
  that
  \begin{multline*}
    \sup_{x\in H} \Bigl\lvert E(f(X^x_T)) - \sum_{(l_1,\ldots,l_p)\in
      \{1,\ldots,N\}^p} \lambda_{l_1} \cdots \lambda_{l_p}
    f(X^x_T(\omega_{l_1,\ldots,l_p})) \Bigr\rvert \\ \leq D T
    \max_{r=1,\ldots,p} (t_r-t_{r-1})^{(m-1)/2}.
  \end{multline*}
\end{proposition}
Note that the assumption on $f$ in Proposition~\ref{prop:bounded_multistep} is
always satisfied if $f$ is bounded, $C^\infty$-bounded and has bounded support
and if all vector fields $\alpha,\beta_1,\ldots,\beta_d $ have bounded support
(compare to Remark \ref{remark_remainder_finite}).  In this case $ P_t f $ has
bounded support, too, and all derivatives are bounded on the whole of $H$ (see
the discussion after Theorem \ref{thm:cubature_dainf} in the next section).
\begin{remark}\label{rem:deterministic_apriori_estimate}
  Proposition~\ref{prop:bounded_multistep} also yields \emph{determinstic a
    priori} estimates for the weak rate of convergence, which hold true if we
  are able to evaluate the respective tree deterministically (see also Section
  \ref{sec:numerical-examples}).  In principle there are methods to do so, see
  \cite{sch/sor/tei:07}, however, they might be more cumbersome to implement
  than a Monte-Carlo evaluation of the tree.
\end{remark}
\begin{remark}
  There is one general case where Proposition \ref{prop:bounded_multistep} can
  be applied: namely for stochastic differential equations of type
  \eqref{eq:sde}. If we replace the unbounded operator $A$ with a bounded
  operator $ \widetilde{A} $, which is close to $ A $ for a large enough set
  of values $ x \in \dom(A) $, then we can apply the previous result for
  bounded linear operators. One candidate for this procedure is the Yosida
  approximation of $A$.
\end{remark}

\section{The case when $A$ is unbounded linear}
\label{sec:cubat-equat-mathc}

We shall also assume in this section that there are no jumps and refer to
equation \eqref{eq:sde}. We recall the analytic background: let us denote by
$\mathcal{D}(A^k)$ the Hilbert space given by $\mathcal{D}(A^k)$ equipped with
the graph norm
\begin{equation*}
  \norm{h}_{\mathcal{D}(A^k)}^2 = \sum_{i=0}^k \norm{A^ih}^2,
\end{equation*}
$h\in\mathcal{D}(A^k)$ and $k\geq 1$. Furthermore, we introduce the space
\begin{equation*}
  \mathcal{D}(A^\infty) = \bigcap_{k=0}^\infty \mathcal{D}(A^k).
\end{equation*}
$\mathcal{D}(A^\infty)$ is topologized as projective limit of the Hilbert
spaces $\mathcal{D}(A^k)$, $k\geq 0$, i.e.~the topology on
$\mathcal{D}(A^\infty)$ is the initial topology of the maps
$\mathcal{D}(A^\infty) \to \mathcal{D}(A^k)$, $k\in \N$. Note that
$\mathcal{D}(A^\infty)$ is no longer a Hilbert space, but only a Fr\'{e}chet
space, i.e.~a locally convex vector space which is completely metrizable by a
translation invariant metric.  We assume our main Assumption
\ref{ass:dainfty}, i.e.~the vector fields restricted to the Sobolev spaces are
$C^\infty$-bounded.

The following proposition collects a few easy, but interesting corollaries
from the existence and uniqueness theorem for equation~(\ref{eq:sde}) applied
to the situation specified in Assumption~\ref{ass:dainfty}.
\begin{proposition}
  \label{prop:solution_dainfty}
  Fix $k\in\N$. For any $x\in \mathcal{D}(A^k)$ there is a unique continuous
  mild solution of~(\ref{eq:sde_nojump}) interpreted as an SDE in the Hilbert
  space $\mathcal{D}(A^k)$. If $x\in \mathcal{D}(A^{k+1}) \subset
  \mathcal{D}(A^k)$, then the mild solution in $\mathcal{D}(A^k)$ coincides
  with the mild solution in $\mathcal{D}(A^{k+1})$. Consequently, for $x\in
  \mathcal{D}(A^{k+1})$, the mild solution of~(\ref{eq:sde}) in
  $\mathcal{D}(A^k)$ is a strong solution and, in particular, a
  semi-martingale.
\end{proposition}

If we start in $x\in \mathcal{D}(A^\infty)$, then we get a continuous process
$X^x_t \in \mathcal{D}(A^\infty)$, such that $X^x_t$ is the (mild and strong)
solution of~(\ref{eq:sde_nojump}) in any Hilbert space $\mathcal{D}(A^k)$,
$k\in\N$. Furthermore, Proposition~\ref{prop:solution_dainfty} allows us to
avoid any problems due to the topological structure of $\mathcal{D}(A^\infty)$
by reverting to an appropriate Hilbert space $\mathcal{D}(A^k)$ and
interpreting the results again in $\mathcal{D}(A^\infty)$. The meaning of $
\dom(A^{\infty}) $ is that it is the largest subspace of the Hilbert space
$H$, where we can innocently do the necessary analysis on differential
operators. Notice that there are subtle phenomena of explosion, which can
occur in this setting: for instance the law of a strong solution process $ X $
solving equation~\eqref{eq:sde} might be bounded in $ H $ but unbounded in $
\dom(A) $, where it is a mild solution. Due to such phenomena, Example
\ref{example_ass:dainf} and Remark \ref{remark_ass:dainf} after Theorem
\ref{thm:cubature_dainf} are in fact quite subtle.

As in Section~\ref{sec:cubat-bound-oper}, we introduce the vector field
$\beta_0$ defined by
\begin{equation}
  \label{eq:strat_corr_unbounded}
  \beta_0(x) = Ax + \alpha(x) - \f{1}{2}\sum_{i=1}^d D\beta_i(x)\cdot
  \beta_i(x). 
\end{equation}
$\beta_0$ is defined for $x\in\mathcal{D}(A)$. As a vector field taking values
in $\mathcal{D}(A^k)$, it is only well-defined on
$\mathcal{D}(A^{k+1})$. Consequently, for $x\in \mathcal{D}(A^{k+1})$, we may
reformulate the SDE~(\ref{eq:sde}) -- understood as equation in
$\mathcal{D}(A^k)$ -- in Stratonovich form~(\ref{eq:sde_stratonovich}).

Now we formulate the stochastic Taylor expansion in some
$\mathcal{D}(A^{r(m)})$ with a degree of regularity $ r(m) \geq 0 $ depending
on $ m \geq 1 $. For the estimation of the error term, we will use the
\emph{extended support} $\esupp(X^x_t;\omega_1,\ldots,\omega_N)$ defined by
\begin{equation}
  \label{eq:extended-support}
  \esupp(X^x_t;\omega_1, \ldots, \omega_N) = \supp(X^x_t) \cup \{
  X^x_t(\omega_1), \ldots, X^x_t(\omega_N)\},  
\end{equation}
where $t > 0$, $x \in H$, and $\omega_1,\ldots,\omega_N$ are paths of bounded
variation. Here $ \supp(X^x_t) $ means the support of the law of $ X^x_t $ in
$\mathcal{D}(A^{r(m)})$.  Despite Assumption~\ref{ass:one-cubature-formula},
let us, for one moment, enter the dependence of the cubature formula on the
interval $[0,t]$ explicitly into the notation, in the sense that
$\omega^{(t)}_1, \ldots, \omega^{(t)}_N$ are the paths of bounded variation
scaled in such a way that they, together with the weights, form a cubature
formula on $[0,t]$. Then we denote
\begin{equation*}
  \s_T(x) = \bigcup_{0\leq s \leq t \leq T} \esupp(X^x_s;
  \omega^{(t)}_1, \ldots, \omega^{(t)}_r).
\end{equation*}
\begin{remark}
  If a general support theorem holds in infinite dimensions, we can replace
  the extended support by the ordinary support of $X^x_t$, since the solution
  of the corresponding ODE driven by paths of bounded variation lie in the
  support of the solution of the SDE according to the support
  theorem. However, up to our knowledge, no general support theorem has been
  established for our setting so far.
\end{remark}

\begin{theorem}\label{thm:stoch_taylor_dainf}
  Let $ m \geq 1 $ be fixed, then there is $ r(m) \geq 0 $ such that for any
  $f \in C^\infty(H; \R)$, $x \in \dom(A^{r(m)})$, and $0<t<1$ we have
  \begin{equation*}
    f(X^x_t) = \sum_{\substack{k\leq m,\, (i_1,\ldots,i_k)\in\mathcal{A}\\
        \deg(i_1,\ldots,i_k) \leq m}} (\beta_{i_1}\cdots \beta_{i_k}f)(x)
    B^{(i_1,\ldots,i_k)}_t + R_m(t,f,x), \quad x \in \mathcal{D}(A^{r(m)}),
  \end{equation*}
  with
  \begin{align*}
    \sqrt{E(R_m(t,f,x)^2)} &\leq C t^{\f{m+1}{2}}
    \max_{m<\deg(i_1,\ldots,i_k)\leq m+2} \; \sup_{0 \leq s \leq t}
    \abs{E(\beta_{i_1}\cdots\beta_{i_k}f(X^x_s))} \\
    &\leq C t^{\f{m+1}{2}} \max_{m<\deg(i_1,\ldots,i_k)\leq m+2} \; \sup_{y
      \in \s_t(x)} \abs{\beta_{i_1}\cdots\beta_{i_k}f(y)}.
  \end{align*}
  We can choose $ r(m) = \floor{\f{m}{2}}+1 $, where $\floor{\f{m}{2}}$ is the
  largest integer smaller than $\f{m}{2}$.
\end{theorem}
\begin{proof}
  The proof is the same as in the finite dimensional situation, but one has to
  switch between different spaces on the way.

  Fix $m$ and $f$ as above and $x \in \mathcal{D}(A^{\floor{\f{m}{2}}+1})$. We
  interpret the equation in $\mathcal{D}(A^{\floor{\f{m}{2}}+1})$. By the
  above remarks, we can express the SDE in its Stratonovich
  form~(\ref{eq:sde_stratonovich}). By It\^{o}'s formula,
  \begin{equation}
    \label{eq:1}
    f(X^x_t) = f(x) + \int_0^t (\beta_0 f)(X^x_s) ds + \sum_{i=1}^d \int_0^t
    (\beta_i f)(X^x_s) \circ dB^i_s.
  \end{equation}
  The idea is to express $(\beta_i f)(X^x_s)$ again by It\^{o}'s formula and
  insert it in equation~(\ref{eq:1}). This is completely unproblematic for
  $i\in \{1,\ldots,d\}$. For $i=0$, recall that
  \begin{equation*}
    (\beta_0 f)(x) = Df(x)\cdot Ax + Df(x)\cdot \Bigl( \alpha(x) -
    \f{1}{2} \sum_{i=1}^d D\beta_i(x)\cdot \beta_i(x)\Bigr).
  \end{equation*}
  By re-expressing~(\ref{eq:1}) in It\^{o} formulation, applying It\^{o}'s
  formula, and re-expressing it back to Stratonovich formulation, we see that
  \begin{equation*}
    (\beta_0 f)(X^x_s) = (\beta_0 f)(x) + \int_0^s (\beta_0^2f)(X^x_u) du +
    \sum_{i=1}^d 
    \int_0^s (\beta_i\beta_0 f)(X^x_u)\circ dB^i_u,
  \end{equation*}
  where
  \begin{equation}
    \label{eq:2}
    (\beta_0^2f)(x) = D^2f(x)(Ax,Ax) + Df(x)\cdot( A^2x + A\alpha(x)+\cdots) +
    \cdots,
  \end{equation}
  provided that all the new vector-fields are well-defined and the processes
  $(\beta_i\beta_0 f)(X^x_u)$ are still semi-martingales. Both conditions are
  satisfied if $x\in \mathcal{D}(A^{2})$ -- notice that the maps
  $\mathcal{D}(A^{k+1}) \to \mathcal{D}(A^k)$, $x\mapsto Ax$ are $C^\infty$,
  $k\in\N$. By induction, we finally get
  \begin{equation*}
    f(X^x_t) = \sum_{\substack{(i_1,\ldots,i_k)\in\mathcal{A}\\
        \deg(i_1,\ldots,i_k) \leq m}} (\beta_{i_1}\cdots \beta_{i_k}f)(x)
    B^{(i_1,\ldots,i_k)}_t + R_m(t,f,x)
  \end{equation*}
  with
  \begin{multline*}
    R_m(t,x,f)\\
    =\sum_{\substack{(i_1,\ldots,i_k)\in\mathcal{A},\,
        i_0\in\{0,\ldots,d\} \\
        \deg(i_1,\ldots,i_k) \leq m < \deg(i_0,i_1,\ldots,i_k)}} \int_{0 \leq
      t_0 \leq \cdots \leq t_k \leq t} (\beta_{i_0}\cdots \beta_{i_k}
    f)(X^x_{t_0}) \circ dB^{i_0}_{t_0} \cdots \circ dB^{i_k}_{t_k}.
  \end{multline*}
  Note that $R_m$ is well-defined for $x \in
  \mathcal{D}(A^{\floor{\f{m}{2}}+1})$ because integration of
  non-semi-martingales with respect to $dt$ is possible, which corresponds to
  the index $i_0=0$.

  As in the finite dimensional case, we re-express $R_m$ in terms of It\^{o}
  integrals and use the (one-dimensional) It\^{o} isometry several times,
  until we arrive at the desired estimate.
\end{proof}

We recall the notation $ P_t f(x) = E(f(X^x_t))$ for bounded measurable
functions $ f: H \to \mathbb{R} $. Analogously to
Proposition~\ref{prop:bounded_multistep}, we immediately get the following
theorem.
\begin{theorem}\label{thm:cubature_dainf}
  Fix $T>0$, $m\geq 1$, $r(m)$, $x \in \dom(A^{r(m)}) $ as in
  Theorem~\ref{thm:stoch_taylor_dainf}, a cubature formula on Wiener space of
  degree $m$ as in Definition~\ref{def:cubature_wiener_space} and a partition
  $0=t_0 < t_1 < \cdots < t_p = T$. Under Assumption~\ref{ass:dainfty}, for
  any $f\in C^\infty(H; \R)$ with
  \begin{equation*}
    \sup_{0\leq t \leq T} \, \sup_{y \in \s_T(x)}
    \abs{\beta_{i_1}\cdots \beta_{i_k} P_tf(y)} < \infty
  \end{equation*}
  for all $(i_1,\ldots,i_k)\in\mathcal{A}$ with $m< \deg(i_1,\ldots,i_k) \leq
  m+2$, $k\in\N$, there is a constant $D$ independent of the partition such
  that
  \begin{multline*}
    \Bigl\lvert E(f(X^x_T)) - \sum_{(l_1,\ldots,l_p) \in \{1,\ldots,N\}^p}
    \lambda_{l_1} \cdots \lambda_{l_p}
    f(X^x_T(\omega_{l_1,\ldots,l_p}))\Bigr\rvert \\
    \leq DT\max_{r=1,\ldots,p} (t_r-t_{r-1})^{(m-1)/2},
  \end{multline*}
  where $X^x_T(\omega)$ is, again, understood as the mild solution to an ODE
  in $\mathcal{D}(A^{r(m)})$ for any path $\omega$ of bounded variation.
\end{theorem}
\begin{proof}
  The proof follows Kusuoka~\cite{kus01}, \cite{kus04}, see also
  \cite{kloe/pla92}. For $f: H\to \R$ and $x \in H$ let
  \begin{equation*}
    Q_{(t)} f(x) = \sum_{l=1}^N \lambda_l f(X^x_t(\omega_l)),
  \end{equation*}
  where $\omega_1, \ldots, \omega_l$ are scaled to form a cubature formula on
  $[0,t]$. Denote $\dt_r = t_r - t_{r-1}$, $r=1, \ldots, p$, the increments of
  the time partition given in the statement of the theorem. By iterating the
  operators $Q_{(\dt_r)}$ (and the semigroup property of ODEs), we immediately
  obtain
  \begin{equation}\label{eq:shorthand-approximation}
    \sum_{(l_1,\ldots,l_p) \in \{1,\ldots,N\}^p} \lambda_{l_1} \cdots \lambda_{l_p}
    f(X^x_T(\omega_{l_1,\ldots,l_p})) = Q_{(\dt_p)} \circ \cdots \circ
    Q_{(\dt_1)} f(x).
  \end{equation}
  By ordinary Taylor expansion, keeping in mind the degree function $\deg$, we
  note that
  \begin{equation*}
    Q_{(t)} f(x) = \sum_{l=1}^N \lambda_l \sum_{\substack{ k\leq m,\,
        (i_1,\ldots,i_k) \in \A \\ \deg(i_1,\ldots,i_k) \leq m}} (\beta_{i_1}
    \cdots \beta_{i_k}f)(x) B^{(i_1,\ldots,i_k)}_t(\omega_l) +
    \widetilde{R}_m(t,x,f), 
  \end{equation*}
  where $x \in \dom(A^{r(m)})$ and
  \begin{multline*}
    \widetilde{R}_m(t,x,f) = \sum_{l=1}^N \lambda_l \sum_{\substack{
        (i_1,\ldots,i_k) \in \A,\, i_0 \in \{0,\ldots,d\}\\
        \deg(i_1,\ldots,i_k) \leq m < \deg(i_0,\ldots,i_k)}}\\ \int_{0 \leq t_0
      \leq \cdots \leq t_k \leq t} (\beta_{i_0} \cdots \beta_{i_k} f)
    (X^x_{t_0}(\omega_l)) d\omega^{i_0}_l(t_0) \cdots d\omega^{i_k}_l(t_k).
  \end{multline*}
  In the sequel $C$ denotes a constant independent of the partition and $x$,
  which may change from line to line. We can estimate the approximation error
  by
  \begin{equation*}
    \abs{\widetilde{R}_m(t,x,f)} \leq C \sup_{0 \leq s \leq t,\, l=1,\ldots, N}
    \; \max_{m \leq \deg(i_0,\ldots,i_k) \leq m+2} \abs{\beta_{i_0} 
      \cdots \beta_{i_k} f(X^x_s(\omega_l))} t^{\f{m+1}{2}}. 
  \end{equation*}
  Combining this result with Theorem~\ref{thm:stoch_taylor_dainf}, we may
  conclude that
  \begin{equation}\label{eq:abs-p_t-fx}
    \abs{ P_t f(x) - Q_{(t)} f(x) } \leq C \sup_{y \in \s_t(x)} \abs{\beta_{i_0}
      \cdots \beta_{i_k} f(y)} t^{\f{m+1}{2}}. 
  \end{equation}
  By telescopic sums,
  \begin{multline*}
    P_T f(x) - Q_{(\dt_p)} \circ \cdots \circ Q_{(\dt_1)} f(x) = \\
    \sum_{r=1}^p Q_{(\dt_p)} \circ \cdots \circ Q_{(\dt_{r+1})} (P_{t_r} f(x)
    - Q_{(\dt_r)} P_{t_{r-1}} f(x)).
  \end{multline*}
  For the estimation of the rear term
  \begin{equation*}
    P_{t_r} f(x) - Q_{(\dt_r)} P_{t_{r-1}} f(x) = (P_{\dt_r} - Q_{(\dt_r)})
    P_{t_{r-1}} f(x),
  \end{equation*}
  we may use~\eqref{eq:abs-p_t-fx} with $f(x)$ being replaced by $P_{t_{r-1}}
  f(x)$, giving us
  \begin{align*}
    \lvert P_T f(x) - &Q_{(\dt_p)} \circ \cdots \circ Q_{(\dt_1)} f(x)\rvert
    \leq
    \sum_{r=1}^p e^{CT} \abs{ P_{t_r} f(x) - Q_{(\dt_r)} P_{t_{r-1}} f(x) } \\
    &\leq C \sum_{r=1}^p \sup_{\substack{y \in \s_{\dt_r}(x)\\ m \leq
        \deg(i_0,\ldots,i_k) \leq m+2}} \abs{\beta_{i_0}
      \cdots \beta_{i_k} P_{t_{r-1}} f(y)} (\dt_r)^{\f{m+1}{2}}\\
    &\leq C \sup_{\substack{y \in \s_T(x),\, 0 \leq t \leq T \\m \leq
        \deg(i_0,\ldots,i_k) \leq m+2 }} \abs{\beta_{i_0} \cdots \beta_{i_k}
      P_t f(y)} \sum_{r=1}^p (\dt_r)^{\f{m+1}{2}},
  \end{align*}
  from which we may easily conclude the theorem.
\end{proof}

\begin{remark}
  Gy\"{o}ngy and Shmatkov~\cite{gyo/shm06} show a strong Wong-Zakai-type
  approximation result, where they also need to impose smoothness assumptions
  on the initial value $x$. Otherwise, the assumptions in~\cite{gyo/shm06} are
  different from ours. They allow linear, densely defined vector fields and
  general adapted coefficients, on the other hand the generator $A$ needs to
  be elliptic.
\end{remark}

\begin{remark}
  Under the previous assumptions we can also prove a Donsker-type result on
  the weak convergence of the ``cubature tree'' to the diffusion. This result
  will be presented elsewhere.
\end{remark}

\begin{remark}
  If $ f $ is smooth then we can show by (first and higher) variation
  processes, as introduced for instance in \cite{for/lue/tei07}, that $ x
  \mapsto P_t f(x) $ is smooth on $ \dom(A^k) $.

  Fix $ k \geq 0 $. Let $J_{0 \to t}(x) \cdot h$ denote the first variation
  process of $X^x_t$ in direction $h \in \dom(A^{k})$, i.e.~
  \begin{equation*}
    J_{0\to t}(x) \cdot h = \left. \f{\pa}{\pa \epsilon} \right|_{\epsilon =
      0} X^{x+\epsilon h}_t \in \dom(A^k).
  \end{equation*}
  $J_{0 \to t}(x) \cdot h$ is the mild solution to an SDE of the
  type~\eqref{eq:sde_nojump}. Consequently, it is bounded in
  $L^2(\Omega,\mathcal{F}, P; \dom(A^k))$ and we may conclude that
  \begin{equation*}
    \left. \f{\pa}{\pa \epsilon} \right|_{\epsilon=0} P_tf(x+\epsilon h) = E
    \Bigl( Df(X^x_t) \cdot J_{0\to t}(x) \cdot h \Bigr)
  \end{equation*}
  exists and is bounded by boundedness of $Df$ and integrability of the first
  variation. Similarly, we get existence and continuity of higher order
  derivatives on $ \dom(A^k) $.
\end{remark}

\begin{example}\label{example_ass:dainf}
  We shall provide examples, where the assumptions of Theorem
  \ref{thm:cubature_dainf} are satisfied, i.e.~where we obtain high-order
  convergence of the respective cubature methods. The conditions seem at first
  sight restrictive (see the following Remark \ref{remark_ass:dainf} for a
  concrete example under Assumptions \ref{ass:dainfty}), however, the
  conditions are parallel to those obtained in \cite{kus04} and
  \cite{lyo/vic04}, where the functions and vector fields have to be bounded
  and $C^{\infty}$-bounded. Here we have an additional complication of one
  certainly unbounded, but not even continuous drift vector field, which leads
  to the following set of assumptions.

  The vector fields $ \alpha, \beta_1,\ldots,\beta_d $ have the following
  property (compare also to tame maps in \cite{Ham82}): there exist smooth
  maps
  \begin{equation*}
    \widetilde{\alpha}, \widetilde{\beta}_1,\ldots,\widetilde{\beta}_d:
    \dom(A^{r(m)}) \to \dom(A^{2r(m)})
  \end{equation*}
  such that for $ k \geq 0 $ the restrictions to the respective subspaces $
  \dom(A^{k+r(m)}) \subset \dom(A^{r(m)}) $ take values in
  $\dom(A^{k+2r(m)})$, i.~e.
  \begin{equation*}
    \widetilde{\alpha}, \widetilde{\beta}_1,\ldots,\widetilde{\beta}_d:
    \dom(A^{k+r(m)}) \to \dom(A^{k+2r(m)}),
  \end{equation*}
  and such that
  \begin{equation}\label{factoring}
    \alpha = \widetilde{\alpha} \circ ({R(\lambda,A)}^{r(m)}),\, \beta_i =
    \widetilde{\beta}_i \circ ({R(\lambda,A)}^{r(m)}), 
  \end{equation}
  for $i=1,\ldots,d$. Here $R(\lambda, A)$ denotes the resolvent map for
  $\lambda \in \rho(A)$. We assume that $
  \widetilde{\alpha},\widetilde{\beta}_1,\ldots,\widetilde{\beta}_d $ have
  bounded support on $ \dom(A^{r(m)}) $. The function $ f $ is of the same
  type $ f = g \circ ({R(\lambda,A)}^{r(m)}) $ for a bounded,
  $C^{\infty}$-bounded function $ g : H \to \mathbb{R} $.

  Under these assumptions we can readily check that the law of the mild
  solution $X^x_t$ starting at the initial value $ x \in \dom(A^{r(m)}) $ has
  bounded support in $ H $: outside a ball of radius $ R > 0 $ in $H$ the
  solution process becomes deterministic, $ X^x_t = S_t x $ on some interval,
  hence by the uniform boundedness theorem there is a large $ R' $ such that
  the image of the ball with radius $ R > 0 $ under the maps $ S_t $ lies in a
  ball with radius $ R' > 0 $ on $ [0,T] $.

  For smooth functions $ f $ of the stated type we then have
  \begin{equation*}
    \sup_{0\leq t \leq T} \sup_{y \in H, \norm{y} \leq R'}
    \abs{\beta_{i_0}\cdots \beta_{i_k} P_tf(y)} < \infty.
  \end{equation*}
  Since we only take the supremum over bounded sets, namely the extended
  supports of $ X^x_t $, this implies the assumption of the Theorem
  \ref{thm:cubature_dainf}.
\end{example}
\begin{remark}\label{remark_ass:dainf}
  The previous assumptions \eqref{factoring} on the vector fields are not too
  restrictive since we can always obtain them by a linear isomorphism and
  (smoothly) cutting off outside a ball in $ \dom(A^{r(m)}) $. Both operations
  are numerically innocent. Under Assumption \ref{ass:dainfty} we can apply
  the following isomorphism to the solution of our SDE \eqref{eq:sde}:
$${R(\lambda,A)}^{-r(m)}: \dom(A^{r(m)}) \to H. $$
This isomorphism transforms the solution $ X^x_t $ on $ \dom(A^{r(m)}) $ to an
$ H $-valued process
$$ Y_t^y = {R(\lambda,A)}^{-r(m)}X^{{R(\lambda,A)}^{r(m)}y}_t $$
with $x = R(\lambda, A)^{r(m)} y$. The transformed process $Y^y_t$ satisfies
an SDE, where the transformed vector fields (if well defined) factor over $
{R(\lambda,A)}^{r(m)} $ such as previously assumed in the assumption
\eqref{factoring}, namely
\begin{multline}
  dY_t^y = (AY^y_{t} + (({R(\lambda,A)}^{-r(m)} \circ \alpha \circ
  {R(\lambda,A)}^{r(m)})(Y^y{t})) dt + \\ + \sum_{i=1}^d
  ({R(\lambda,A)}^{-r(m)} \circ \beta_i \circ {R(\lambda,A)}^{r(m)})(Y^y_{t})
  dB^i_t.
\end{multline}
The assumptions \eqref{factoring} mean that we must (smoothly) cut off the
vector fields $ \alpha, \beta_1,\ldots,\beta_d $ outside sets of large norm
${||.||}_{\dom(A^r(m))}$, which is an event -- under Assumption
\ref{ass:dainfty} -- of small probability (recall that the vector fields are
Lipschitz on $ \dom(A^{r(m)}) $ and therefore second moments with respect to
the norm ${||.||}_{\dom(A^r(m))}$ exist). Notice that the cut-off vector
fields do not have an extension to $ H $ since continuous functions with
bounded support on $ \dom(A^{r(m)}) $ do generically not have a continuous
extension on $ H $. For $ {(Y_t^y)}_{t \geq 0} $ we can take initial values $
y \in H $, however, those initial values correspond to quite regular initial
values $ x = {R(\lambda,A)}^{r(m)} y \in \dom(A^{r(m)}) $ for the original
process $ {(X_t^x)}_{t \geq 0} $.

From the point of view of the process $ Y $ we have hence proved that
$$
f (Y_T^y) = g \circ ({R(\lambda,A)}^{r(m)}) \circ ({R(\lambda,A)}^{-r(m)})
(X_T^x) = g (X^x_T)
$$
is weakly approximated by evaluating $ f $ on the cubature tree for $ Y
$. This is equivalent to evaluating $ g $ on the cubature tree for $ X $ in
order to approximate the expected value $ E(g(X^x_T))$.
\end{remark}

\section{The cubature method in the presence of jumps}
\label{sec:jump-structure}
The extension of cubature formulas to jump diffusions seems to be new even in
the finite dimensional case. We shall heavily use the fact that only finitely
many jumps occur in compact time intervals almost surely.

We shall first prove an asymptotic result on jump-diffusions.

\begin{proposition}\label{conditioning on jumps}
  Consider equation \eqref{eq:sde-j}. Let $ f: H \to \mathbb{R} $ be a bounded
  measurable function, then we obtain
  \begin{multline}
    E(f(X^x_t)) = \sum_{n_{1},\dots,n_{e} \geq 0} \frac{_{\mu_{1}^{n_{1}}
        \cdots \mu_{e}^{n_{e}}}}{n_{1}!\cdots n_{e}!}e^{- t \mu_{1} n_{1}-
      \ldots - t \mu_{e} n_{e}}
    t^{n_1 + \ldots + n_e} \times \\
    \times E(f(X^x_t) \; |\; N_{t}^{j}=n_{j}\text{ for }j=1,\dots,e)
  \end{multline}
  for $ t \geq 0 $.
\end{proposition}

\begin{proof} 
  We condition on the jump times and read off the results by inserting the
  probabilities for a Poisson process with intensity $ \mu_j $ to reach level
  $ n_j $ at time $ t $.
\end{proof}

This result gives us the time-asymptotics with respect to the
jump-structure. It can now be combined with the original cubature result for
the diffusion between the jumps, in order to obtain a result for
jump-diffusions.  We denote by $ \tau^j_n $ the jump-time of the Poisson
process $ N^j $ for the $n$-th jump. We know that for each Poisson process the
vector $ (\tau^j_1,\ldots,\tau^j_{n_j} - \tau^j_{{n_j}-1},t - \tau^j_{n_j}) $
is uniformly distributed if conditioned on the event that $ N^j_t = n_j \geq 1
$. The uniform distribution is on the $n_j$-simplex $t\Delta^{n_j} \subset
\mathbb{R}^{n_j+1} $. This allows us to apply an original cubature formula
between two jumps of order $ m - 2n_1 - \ldots -2n_e$, since we gain, for each
jump, one order of time-asymptotics from the jump structure.

Assume now that the jump distributions $ \nu_j $ are concentrated at one point
$ z_j \neq 0 $, i.e.~$ \Delta L_{\tau^j_k}^j = z_j $ for $ j=1,\ldots,e$ and $
k \geq 1$. If we want to consider a general jump-structure this amounts to an
additional integration with respect to the jump distribution $ \nu_j $.

Now we define an short-time approximation for the conditional expectations
\begin{equation*}
  E(f(X^x_t) \; |\; N_{t}^{j}=n_{j}\text{ for }j=1,\dots,e)
\end{equation*}
of order $ m - 2n_1 - \ldots - 2n_e $ with $ n_1 + \ldots + n_e \leq
\frac{m+1}{2} $. Expressed in words, we are going to do the following:
starting from the initial value $ x \in \dom(A^{r(m)}) $ we solve the
stochastic differential equation~\eqref{eq:sde-j} along the cubature paths $
\omega_l $ with probability $ \lambda_l > 0 $, $l=1,\ldots,N$ until the first
jump appears. We collect the end-points of the trajectories, add the jump size
at these points and start a new cubature method from the resulting points
on. Notice that we can take a cubature method of considerably lower degree
since every jump increases the local order of time-asymptotics by $ 1 $. The
jump times are chosen independent and uniformly distributed on simplices of
certain dimension $ n_j $ such that $ n_1 + \ldots + n_e = n \leq
\frac{m+1}{2} $. We denote the cubature trajectory between jump $ \tau^j_{q-1}
$ and $ \tau^j_q $ for $ 1 \leq q \leq n_j $ with $ \omega_{l,j,q} $. If $ n_j
= 0 $ no trajectories are associated. Hence we obtain the following theorem:
\begin{theorem}
  Fix $ m \geq 1 $. Consider the stochastic differential equation
  (\ref{eq:sde-j}) under the condition $ N^j_t = n_j $ for $ j=1,\ldots,m $
  with $ n_1 + \ldots + n_e = n $ along concatenated trajectories of type $
  \omega_{l,j,q} $. Choose a cubature method of degree
  \begin{equation*}
    m' = m - 2n \geq 1.
  \end{equation*}
  Concatenation is only performed with increasing $ q $-index and a typical
  concatenated trajectory is denoted by $ \omega_{l_1,\ldots,l_n} $. Here we
  have in mind that the intervals, where the chosen path is $ \omega_{l,j,q}
  $, come from a jump of $N^j$ and have length $ \tau_{q}^j - \tau^j_{q-1} $.
  Then there is $ r(m') \geq 0 $ such that
  \begin{multline}
    \bigl| E(f(X^x_t)\; |\; N_{t}^{j}=n_{j}) \bigr. \bigl. -
    \sum_{l_1,\ldots,l_n =1}^N \lambda_{l_1} \dots \lambda_{l_n}
    E(f(X^x_t(\omega_{l_1,\ldots,l_n}))\; |\; N_{t}^{j}=n_{j}) \bigr| \\ 
    \leq C t^{\f{m'+1}{2}} \max_{\substack{(i_1,\ldots,i_k)\in \mathcal{A} \\
        \deg(i_1,\ldots,i_k)\leq m+2}} \sup_{\substack{y \in
        \operatorname{supp}(X^x_s), \\ 0 \leq s \leq t}}
    |\beta_{i_1} \cdots \beta_{i_k} E(f(X^y_{{\tau_q},t} (
    \omega_{l_{q+1},\ldots,l_{n}} ))\,  
    | \, N_t^j = n_j)|,
  \end{multline}
  where $ X^x_t(\omega) $ means the solution of the stochastic differential
  equation~\eqref{eq:sde-j} in Strato\-novich form
  \begin{equation}
    dX_t^x(\omega) = \beta_0(X_{t^-}^x(\omega))dt + \sum_{i=1}^d
    \beta_i(X^x_{t^-}(\omega)) \circ dB^i_t +
    \sum_{j=1}^{e}\delta_{j}(X_{t^-}^{x})dL_{t}^{j}, 
  \end{equation}
  along the trajectory $ \omega $.
\end{theorem}
\begin{proof}
  By our main Assumption \ref{ass:dainfty} we know that the linkage operators
  $ x \mapsto \delta^j(x) $ are $ C^{\infty}$-bounded on each $ \dom(A^k) $,
  hence through concatenation the errors, which appear on each subinterval $
  [t_{n-q},t[ $ are of the type $ y \mapsto
  E(X^y_{\tau_{q},t}(\omega_{l_{q+1}q,\ldots,l_{n}})) $ for some $ 1 \leq q
  \leq n $. Taking the supremum yields the result.
\end{proof}

Combining the previous result with Proposition \ref{conditioning on jumps}
yields under certain conditions on the vector fields (see the discussion after
Theorem \ref{thm:stoch_taylor_dainf} in the previous section) by the triangle
inequality that there is a constant $ D >0 $ such that
\begin{multline}
  \bigl| E(f(X^x_t)) - \sum_{2(n_1+\ldots+n_e) \leq m} \sum_{l_1,\ldots,l_n
    =1}^N \frac{_{\mu_{1}^{n_{1}} \cdots \mu_{m}^{n_{e}}}}{n_{1}!\cdots
    n_{e}!}e^{- t \mu_{1} n_{1}- \ldots - t \mu_{e} n_{e}} \times \bigr. \\
  \bigl. \times \lambda_{l_1} \dots \lambda_{l_n}
  E(f(X^x_t(\omega_{l_1,\ldots,l_n}))\; |\; N_{t}^{j}=n_{j}) \bigr| \leq D
  t^{\f{m+1}{2}} < \infty,
\end{multline}

By iteration of the previous result we obtain in precisely the same manner as
in Section \ref{sec:cubat-equat-mathc} a cubature method of order $ m $ by
applying several cubature methods of order $ m' \leq m $ between the jumps.

\begin{remark}
  The only random element in the expectation
  $E(f(X^x_t(\omega_{l_1,\ldots,l_n}))\; |\; N_{t}^{j}=n_{j})$ is given by the
  jump times $ \tau^j_k $, which vary on certain simplices. For the
  implementation one has to simulate the uniform distributions on the
  simplices $ t\Delta^k $. Since the integrals on the simplices $ t\Delta^k $
  only have continuous integrands, we cannot hope for other methods than
  Monte-Carlo. The evaluation by a Monte-Carlo-algorithm can also be seen as a
  random choice of the concatenation grid for the constructed
  $\omega_{l_1,\ldots,l_n}$. Another view could be to see a deterministic grid
  for the diffusion which is saturated by points where jumps
  occur. Implementation will be done elsewhere.
\end{remark}

\section{Numerical examples}
\label{sec:numerical-examples}

Due to the previous results a numerical scheme for equation \eqref{eq:sde} can
be set up by the following steps. Notice that we have a weak order of
approximation $\frac{m-1}{2} $ only under the assumptions of the previous
section. In order to obtain those assumptions one has to modify a general
equation of type \eqref{eq:sde} by smoothing procedures. These modifications
can be done in a controlled way, more precisely, for each modification we have
a rate of convergence to the un-modified object.
\begin{itemize}
\item Approximate the vector fields $ \alpha, \beta_1,\ldots,\beta_d $ by
  vector fields satisfying the Assumptions \ref{ass:dainfty}. If the original
  vector fields are globally Lipschitz one can do this approximation with a
  rate of convergence for the $ L^2 $-distance of the original solution
  process and its approximation.
\item Choose a degree of accuracy $ m \geq 2 $, which determines the weak
  order of convergence $ \frac{m-1}{2} $ in the sequel. Associated to $ m $
  the number $ r(m) $ can be identified, which tells us about the degree of
  regularity of the initial value $ x \in H $, which one needs for the
  assertions of Theorem \ref{thm:cubature_dainf}.
\item Identify due to the previous specifications a radius $ R > 0 $ such that
  the $ \dom(A^{r(m)})$-norm of $ X^x_T $ is rarely beyond $ R $. Cut-off the
  vector fields smoothly on $ \dom(A^{r(m)}) $ and verify assumptions
  \eqref{factoring} -- maybe after smoothing -- for the transformed process $
  Y $ such as exercised in Remark \ref{remark_ass:dainf}.
\item The resulting tree of trajectories yields a finite number of
  non-autonomous PDEs \eqref{eq:cubature_ode}, which have to be evaluated. In
  the implementation one calculates with the smoothened vector fields
  satisfying the Assumptions \ref{ass:dainfty}, but not with the cut-off
  vector fields, since a large $ \dom(A^{r(m)})$-norm value $R$ is reached
  with small probability due to its very choice.
\end{itemize}

We test the cubature method for two concrete examples: one toy example, where
explicit solutions of the SPDE are readily available, and another, more
interesting but still very easy example. Since cubature on Wiener space is a
weak method, we calculate the expected value of a functional of the solution
to the SPDE in both cases, i.e.~the outputs of our computations are real
numbers.

The results presented here are calculated in MATLAB using the built-in
PDE-solver \verb|pdepe| for solving the deterministic PDEs given by inserting
the cubature paths into the SPDE under consideration. This PDE solver depends
on a space grid given by the user as well as on a time grid, which is not very
critical because it is adaptively refined by the program.

We do not use recombination techniques for cubature on Wiener space as
in~\cite{sch/sor/tei:07} and use the simplest possible cubature formula for
$d=1$ Brownian motions:
\begin{equation*}
  \omega^{(T)}_1(t) = - \f{t}{\sqrt{T}}, \ \omega_2^{(T)} = \f{t}{\sqrt{T}},\quad
  t\in [0,T],
\end{equation*}
with weights $\lambda_1 = \lambda_2 = \f{1}{2}$ define a cubature formula of
degree $m = 3$ on $[0,T]$. Consequently, solving an SDE on a Hilbert space
with cubature on Wiener space for the above cubature formula and $p$
iterations means solving $2^p$ PDEs. This starts to get restrictive even for a
very simple problem for, say, $p = 10$ -- where one already has to solve more
than one thousand PDEs. One possibility to overcome these tight limitations is
to use ``a Monte-Carlo simulation on the tree''. Recall that an $p$-step
cubature method approximates
\begin{equation*}
  E(f(X^x_T)) \approx \sum_{(j_1,\ldots,j_p) \in \{1,\ldots,N\}^p} \lambda_{j_1}
  \cdots \lambda_{j_p} f(X^x_T(\omega_{j_1,\ldots,j_p})).
\end{equation*}
Since $\sum \lambda_{j_1}\cdots \lambda_{j_p} = 1$, we can interpret the right
hand side as the expectation of a random variable $f(X^x_T(\omega_{\cdot}))$
on the tree $\{1,\ldots,N\}^p$. Therefore, we can approximate the right hand
side by picking tuples $(j_1,\ldots,j_p) \in \{1,\ldots,N\}^p$ at random --
according to their probabilities $\lambda_{j_1}\cdots \lambda_{j_p}$ -- and
calculating the average of the corresponding outcomes
$f(X^x_T(\omega_{j_1,\ldots,j_p}))$. Of course, by following this strategy we
have to replace the deterministic error estimates by a stochastic one, which
heavily depends on the standard deviation of $f(X^x_T(\omega_\cdot))$
understood as a random variable on the tree. Notice, however, that this is the
usual situation for weak approximation methods for SDEs.
 
Consider the Ornstein-Uhlenbeck process $X^x_t$ defined as solution to the
equation
\begin{equation}
  \label{eq:example_ou}
  dX^x_t = \Delta X^x_t dt + \phi dB_t
\end{equation}
on the Hilbert space $H = L^2(]0,1[)$. $\Delta$ denotes the Dirichlet
Laplacian on $]0,1[$, i.e.~$\Delta$ is a negative definite self-adjoint
operator on $H$ with $\mathcal{D}(\Delta) = H^1_0(]0,1[) \cap H^2(]0,1[)$
extending the classical Laplace operator defined on $C^\infty_c(]0,1[)$. It is
easy to see that $\Delta$ is dissipative and therefore, by the Lumer-Phillips
theorem, it is the generator of a $C_0$ contraction semigroup $(S_t)_{t\geq0}$
on $H$. The coefficient $\phi\in H$ is some fixed vector.

In this case, the definition of a mild solution
\begin{equation}
  \label{eq:5}
  X_t^x = S_t x + \int_0^t S_{t-s} \phi dB_s
\end{equation}
already gives a representation of the solution provided that the
heat-semigroup $S_t$ applied to the starting vector $x$ and to $\phi$ is
available. We choose $x(u) = \sin(\pi u)$, $u \in]0,1[$, and may conclude that
\begin{equation*}
  S_t x = e^{-\pi^2 t} x
\end{equation*}
because $x$ is an eigenvector of $\Delta$ with eigenvalue $-\pi^2$. Consider
the linear functional $\Phi : H \to \R$ given by
\begin{equation}
  \label{eq:4}
  \Phi(y) = \int_0^1 y(u) du, \quad y \in H.
\end{equation}
We want to compute
\begin{align*}
  E(\Phi(X^x_1)) &= E\Bigl( \int_0^1 e^{-\pi^2}\sin(\pi u) du +
  \int_0^1 \int_0^1 S_{1-s}\phi(u)dB_s du \Bigr)\\
  &= \int_0^1 e^{-\pi^2} \sin(\pi u) du = 0.3293 \times 10^{-4}.
\end{align*}
\begin{table}[htbp]
  \centering
  \begin{tabular}{|c|c|}
    \hline
    $p$ & Error\\
    \hline
    $1$ & $-0.3601 \times 10^{-4}$\\
    $2$ & $-0.2192 \times 10^{-4}$\\
    $3$ & $-0.1226 \times 10^{-4}$\\
    $4$ & $-0.0652 \times 10^{-4}$\\
    $5$ & $-0.0334 \times 10^{-4}$\\
    $6$ & $-0.0172 \times 10^{-4}$\\
    $7$ & $-0.0084 \times 10^{-4}$\\
    $8$ & $-0.0031 \times 10^{-4}$\\
    $9$ & $-0.0002 \times 10^{-4}$\\
    $10$ & $-0.0013 \times 10^{-4}$\\
    \hline
  \end{tabular}
  \caption{Error for the cubature method in the OU-case (absolute error)}
  \label{tab:cub50vas}
\end{table}
\begin{figure}[htbp]
  \centering
  \includegraphics[width=\textwidth]{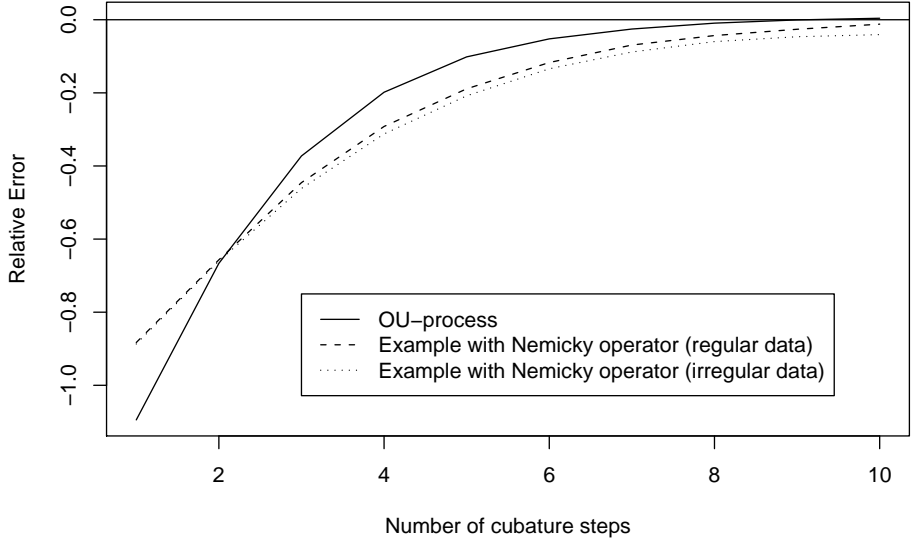}
  \caption{Relative errors for the OU-process (equation~(\ref{eq:example_ou}))
    and the process with nonlinear volatility
    (equation~(\ref{eq:heat_nemicky})). In the latter case, both initial
    values $x(u) = \sin(\pi u)$ -- referred to as ``regular case'' -- and $x$
    given by~(\ref{eq:6}) -- referred to as ``irregular case'' -- are used.}
  \label{fig:cub50vas}
\end{figure}

In Table~\ref{tab:cub50vas}, the error, i.e.~the output of the method minus
the true value given above, is presented. $p$ is the number of cubature steps,
i.e.~the number of iterations of the one-step cubature method. The
discretization in space, i.e.~of $]0,1[$, used by the PDE solver contains $50$
uniform points, the discretization of the time-interval -- additional to the
one induced by the cubature method -- contains $500$ points.  The stochastic
perturbation factor $\phi$ is chosen to be $\phi(u) = \sin(\pi u)$, i.e.~$\phi
\in \mathcal{D}(\Delta^\infty)$ even. We see a very fast decrease of the error
in this simple situation. On the other hand, the variance of the random
variable on the tree considered before is clearly too high for the Monte-Carlo
simulation on the tree to work. Indeed, $\Phi(X^x_1)$ has true standard
deviation of
\begin{equation}
  \label{eq:ou-sd}
  \operatorname{sd}(\Phi(X^x_1)) = \sqrt{\frac{2}{\pi^4} ( 1 - e^{-2\pi^2})} =
  0.1433. 
\end{equation}
Assuming that the central limit theorem applies, confidence intervals around
the solution given by a Monte-Carlo method are proportional to the standard
deviation divided by the square root of the number of
trajectories. Consequently, we would roughly need to calculate $10^{12}$ paths
on the tree in order to achieve a similar level of exactness as in
Table~\ref{tab:cub50vas}! Indeed, note that the standard deviation of the
solution is of order $10^{-1}$, while the error in the last row of
Table~\ref{tab:cub50vas} is of order $10^{-7}$. The equation
\begin{equation*}
  \frac{10^{-1}}{\sqrt{m}} \approx 10^{-7}
\end{equation*}
then gives $m \approx 10^{12}$. Note that this heuristics is also confirmed by
our experiments, where Monte-Carlo simulation on the tree clearly fails. The
data are also shown in Figure~\ref{fig:cub50vas}.

\begin{remark}
  \label{mc-dooes-not-work}
  The failure of Monte-Carlo simulation on the tree also applies to any other
  (naive) Monte-Carlo approach to problem~(\ref{eq:example_ou}), including the
  usual finite element or finite difference approaches.
\end{remark}

As a more realistic example we consider the heat equation with a stochastic
perturbation involving a Nemicky operator. More precisely, consider
\begin{equation}\label{eq:heat_nemicky}
  dX^x_t = \Delta X^x_t dt + (\sin \circ X^x_t) dB_t,
\end{equation}
with $x(u) = \sin(\pi u)$. Even though we do not know the law of the solution
$X^x_1$ of~(\ref{eq:heat_nemicky}), we are still able to calculate
$E(\Phi(X^x_1))$ explicitly because $\Phi$ is a linear functional. Indeed,
$X^x_1$ is given by
\begin{equation}
  X^x_1 = S_1x + \int_0^1 S_{1-s} (\sin \circ X^x_s) dB_s
\end{equation}
and, consequently,
\begin{equation*}
  \Phi(X^x_1) = \Phi(S_1x) + \int_0^1 \Phi(S_{1-s} \sin \circ X^x_s) dB_s.
\end{equation*}
The expectation of the (one-dimensional) It\^{o}-integral is $0$ and we get
the same result as before, i.e.
\begin{equation*}
  E(\Phi(X^x_1)) = \Phi(S_1x) = 0.3293\times 10^{-4}
\end{equation*}
for $x(u) = \sin(\pi u)$. Nevertheless, we believe that this example is
already quite difficult, especially since the cubature method actually has to
work with the Stratonovich formulation
\begin{equation}
  \label{eq:heat_nemicky_strat}
  dX^x_t(\omega) = \bigl(\Delta X^x_t(\omega) - \f{1}{2} (\cos \circ
  X^x_t(\omega)) (\sin\circ 
  X^x_t(\omega))\bigr) dt + (\sin \circ X^x_t(\omega)) d\omega(t).
\end{equation}
In particular, the equation (in Stratonovich form) has a non-linear drift and
a non-linear volatility.

Note that we expect the standard deviation of the solution of the above
equation to be smaller than before, because $(\sin \circ X^x_t)^2$ decreases
as $X^x_t$ decreases in $t$.
\begin{table}[htbp]
  \begin{tabular}{|c|c|}
    \hline
    $l$ & Error \\
    \hline
    $1$ & $-0.2907\times 10^{-4}$\\
    $2$ & $-0.2163\times 10^{-4}$\\
    $3$ & $-0.1467\times 10^{-4}$\\
    $4$ & $-0.0961\times 10^{-4}$\\
    $5$ & $-0.0622\times 10^{-4}$\\
    $6$ & $-0.0385\times 10^{-4}$\\
    $7$ & $-0.0228\times 10^{-4}$\\
    $8$ & $-0.0142\times 10^{-4}$\\
    $9$ & $-0.0086\times 10^{-4}$\\
    $10$ & $-0.0040\times 10^{-4}$\\
    \hline
  \end{tabular} 
  \caption{Results of the cubature method for~(\ref{eq:heat_nemicky}) (absolute
    error)} 
  \label{tab:cub_nem}
\end{table} 
The space discretization used by the PDE-solver has size $50$, which already
seems to be sufficient, because using a finer discretization ($100$ grid
points) does not change the results significantly.
\begin{table}[htbp]
  \begin{tabular}{|c|c|c|c|}
    \hline
    $l$ & $m$ & Error & Stat.~Error \\
    \hline
    $ 5$ & $  32$ & $\phantom{-}0.0567\times 10^{-4}$ & $0.1498\times 10^{-4}$ \\
    $10$ & $1000$ & $-0.0325\times 10^{-4}$ & $0.0179\times 10^{-4}$ \\
    $15$ & $1500$ & $-0.0184\times 10^{-4}$ & $0.0172\times 10^{-4}$ \\
    $20$ & $2000$ & $\phantom{-}0.0128\times 10^{-4}$ & $0.0170\times 10^{-4}$ \\
    $25$ & $2500$ & $\phantom{-}0.0179\times 10^{-4}$ & $0.0145\times 10^{-4}$ \\
    $30$ & $3000$ & $\phantom{-}0.0596\times 10^{-4}$ & $0.0167\times 10^{-4}$ \\
    \hline
  \end{tabular} 
  \caption{Results of the cubature method with Monte-Carlo simulation on the 
    tree for~(\ref{eq:heat_nemicky}) (absolute error)}
  \label{tab:mccub_nem}
\end{table} 
Table~\ref{tab:mccub_nem} shows the results using Monte-Carlo simulation on
the tree. $m$ denotes the number of trajectories followed, while the
``Statistical Error'' in the table is an indicator for the error of the
Monte-Carlo simulation. More precisely, the values in the last column are the
empirical standard deviations of the result divided by the square root of the
number of trajectories. Comparable to the Ornstein-Uhlenbeck process, the
convergence of the pure cubature method is very fast, see
Table~\ref{tab:cub_nem}. The (empirical) variance is, however, quite large
such that the Monte-Carlo aided method does not work at all. Note that the
statistical error in Table~\ref{tab:mccub_nem} is of the order of the total
computational error, which can be almost completely attributed to the Monte
Carlo simulation.

To test the method further we also try more irregular data. Let
\begin{equation}
  \label{eq:6}
  x(u) = \f{1}{2} \sqrt{\f{1-2\abs{u-\f{1}{2}}}{\sqrt{\abs{u-\f{1}{2}}}}}.
\end{equation}
The exact value of the quantity of interest $E(\Phi(X^x_1)) = \Phi(S_1x)$ is
calculated by solving the corresponding heat equation numerically. This gives
the value $E(\Phi(X^x_1)) = 0.3002\times 10^{-4}$.  The initial vector $x$
given in~(\ref{eq:6}) is in $L^2(]0,1[)$ but its derivative is no longer
square-integrable. Consequently, $x\notin \mathcal{D}(A)$ and the theory does
not provide an order of approximation. Nevertheless, probably due to the
smoothing-properties of the Laplace operator numerical results show the same
behavior as before, see Figure~\ref{fig:cub50vas}.

If we replace the heat equation~(\ref{eq:example_ou}) by an evolution equation
of the form
\begin{equation}
  \label{eq:7}
  dX^x_t = \f{d}{du} X^x_t dt + \sin\circ X^x_t dB_t,
\end{equation}
then we still see the same kind of behavior if we fix the space-discretization
for the PDE-solver. This time, the PDEs require a much finer space resolution
in order to give reliable numbers.

\end{document}